\documentclass[10pt]{article}
\usepackage{amsmath, dsfont}
\usepackage{amssymb}
\usepackage[dvips]{graphics}
\usepackage{epsfig}
\usepackage{calc}
\usepackage{amsfonts}
\usepackage{graphicx}
\usepackage[spanish]{babel}
\usepackage[usenames,dvipsnames]{pstricks}
\usepackage[ansinew]{inputenc}

\title{La derivada en varias variables como analog\'ia formal de su hom\'onima escalar\\
}

\author{M. Bravo-Gaete, F. Córdova-Lepe, P. Dotte. \\
{\it Facultad de Ciencias B\'asicas, Universidad Cat\'olica del
Maule,}\\ {\it 3605 San Miguel Av., Talca, Chile}}

\begin{document}
\maketitle

\begin{abstract}
Algunas {\it analog\'ias formales} entre el C\'alculo Diferencial en
Una Variable y el C\'alculo Diferencial en Varias Variables son
presentadas. Se estudia e introduce la {\it derivabilidad de
funciones de varias variables} desde su an\'alogo conceptual
escalar. Lo anterior, explorando la imagen din\'amica de l\'imite de
una familia de pendientes de planos secantes al gr\'afico de una
funci\'on bivariada.

\vspace{2mm}

\end{abstract}

\section{Introducci\'on}
Los procesos de ense\~nanza, aprendizaje y comunicaci\'on de la
matem\'atica propia del modelamiento de la variaci\'on y el cambio
(entendidos en contextos de un sistema educativo y un medio social)
conforman la materia de estudio del  {\it Pensamiento Variacional}
(PV). Al entender el {\it cambio} en un sistema como la
constataci\'on del paso de un estado a otro y la {\it variaci\'on}
como la cuantificaci\'on del tama\~no de dicho paso, es claro que un
tipo de objetos claves del PV son los {\it operadores de
diferenciaci\'on}, en particular y privilegiadamente, para cierta
clase de modelos, la derivada.

A todo nivel educativo, la tarea de construir e instalar en
profundidad en los estudiantes el concepto de variaci\'on y sus
respectivas lecturas matem\'aticas se considera bastante dif\'icil,
G\'OMEZ \cite{GomezO}.  En general, esta labor requiere la
articulaci\'on e integraci\'on de muy distintos saberes y modos de
registro. Aunque el caso de la derivada de funciones (reales de una
variable real) en este sentido est\'a bastante documentado, ARTIGUE
\cite{Artigue}, no ocurre lo mismo al pasar a varias variables,
espec\'ificamente con la noci\'on de diferenciabilidad.

Situados en la educaci\'on terciaria, existe un largo camino de
previos conceptuales, que van desde el Axioma del Supremo hasta el
L\'imite de Funciones, que deber\'ian estar consolidados antes de
arribar en el desarrollo de los cursos a la noci\'on de la derivada.
Sin este encadenamiento de conceptos, el generar buenos niveles de
comprensi\'on y de competencias dirigidas a la construcci\'on y
an\'alisis de algunos tipos de modelos (v.g., deterministas y
continuos de sistemas de naturaleza din\'amica) es un quehacer muy
complejo, ver S\'ANCHEZ-MATAMOROS y col. \cite{Sanchez}. La
modelizaci\'on tendiente a clarificar el entendimiento y a
introducir predicci\'on de los cambios de estado, pasa previamente
por un constructo, la instalaci\'on del pensamiento y el lenguaje
variacional vinculado a la derivada. En cuanto a debilidades
formativas en esta direcci\'on, DOLORES en \cite{Dolores00} expresa:

\begin{center}
{\it {``...Los investigadores en este campo coinciden en que,
cantidades significativas de estudiantes s\'olo pueden obtener
derivadas de funciones algebraicas mediante f\'ormulas, pero
dif\'icilmente comprenden el para qu\'e de esos algoritmos que
realizan y el significado de los conceptos. Inclusive,
dif\'icilmente logran asociar las ideas claves del c\'alculo en la
resoluci\'on de problemas elementales sobre la variaci\'on, a pesar
de que hist\'oricamente del estudio de estos \'ultimos se originaron
las ideas claves del CD''.}}
\end{center}

Nuestro trabajo se enmarca en el PV en cuanto exploramos y
proponemos uno modo de construcci\'on de objetos matem\'aticos
asociados al tr\'ansito desde el C\'alculo Diferencial en Una
Variable (CALUNO) hacia el C\'alculo Diferencial en Varias Variables
(CALVAV). M\'as precisamente, nos referimos a una alternativa, con
indicios pedag\'ogicos, \'util para una presentaci\'on dial\'ogica
del proceso de generalizar a funciones bivariadas los conceptos de
{\it derivada} y {\it diferenciabilidad}. Lo anterior, motivados por
la fuerte presencia del C\'alculo Diferencial (CD) en las
curr\'icula formativas de gran parte de los programas de formaci\'on
profesional en las instituciones de ense\~nanza superior y, muy
especialmente, en intentar destrabar las consabidas dificultades que
representa para los estudiantes el aprendizaje del C\'alculo.

Es muy com\'un que la ense\~nanza del CD sea llevada en forma
abrupta y no relacionada con las tem\'aticas y conocimientos previos
que podr\'ian cumplir el rol de puente, ver DOLORES
\cite{Dolores96}. Adem\'as, en el CD a veces, el \'enfasis se
concentra en los procedimientos algebraicos, sin llegar a evidenciar
la din\'amica del cambio que involucra el pensamiento variacional.
El tratamiento en varias variables no est\'a ajeno a estos problemas
y el paso de CALUNO a CALVAV puede llegar a convertirse en un enorme
salto conceptual. De lo que es la experiencia con la derivada,
pasamos a tipos de derivadas (parciales o direccionales) y de ah\'i
hacia otro abismo, la idea de diferenciabilidad.

En los primeros cursos de CD (i.e., en una variable) el internalizar
la derivada consume buena parte del tiempo de los  maestros. La
implementaci\'on de iniciativas did\'acticas que logren un
afiatamiento y robustez m\'inima del concepto en sus estudiantes
requiere tiempo y esfuerzo. Aunque, despu\'es de esta etapa los
profesores solemos operacionalizar mostrando propiedades para el
c\'alculo expl\'icito de derivadas y, finalmente, con frecuencia nos
centramos en sus aplicaciones. La primera etapa del proceso anterior
se reconoce exitosa cuando el estudiante logra describir (verbal,
gr\'afica y simb\'olicamente) la derivada de una funci\'on en un
punto como el l\'imite (de existir) de una familia de pendientes de
rectas secantes sujetas (que contienen) a dicho punto. Lo que en
particular, permite introducir la formalizaci\'on de la idea de
recta tangente al gr\'afico de la funci\'on en un punto dado. Todo
este proceso, algo m\'as o un poco menos, se tiende a reproducir
cuando se enfrentan a la definici\'on de la derivada direccional en
un primer curso de CALVAV. En general, para una funci\'on bivariada
se espera la visualizaci\'on, ahora en tres dimensiones, de rectas
secantes a la gr\'afica de una curva que es la intersecci\'on entre
la gr\'afica de la funci\'on en cuesti\'on y el plano vertical al
dominio definido por el vector direcci\'on.

En este trabajo, proponemos la v\'ia de la construcci\'on de
semejanzas con soporte en los conocimientos previos de la derivada
escalar. En este caso, se usan los conocimientos de CALUNO como el
dominio base de {\it analog\'ias} que transportan razgos y
propiedades hacia lo que debiera ser la noci\'on de derivada de
funciones bivariadas y, por simple generalizaci\'on, a funciones
reales de m\'as de dos variables.

M\'as precisamente, pretendemos que gran parte de la experiencia de
los estudiantes sobre la instalaci\'on en ellos de la derivada de
funciones reales univariadas persista mediante una ``buena''
analog\'ia al subir dimensionalmente al caso de funciones de dos
variables. De modo tal que, la validez sint\'antica de las
expresiones matem\'aticas y la sem\'antica de los relatos y eslogans
persistan y puedan ayudarnos en la aventura de adentrarnos, junto al
estudiante, en el nuevo territorio, el de los operadores de
derivaci\'on de las funciones bivariadas. Como hemos dicho, la
t\'ecnica ser\'a la construcci\'on de analog\'ias, pero en el
sentido m\'as riguroso, introducido en C\'ORDOVA-LEPE et al.
\cite{Cordova2}, y denominadas {\it analog\'ias formales}, esto es,
una sistematizaci\'on y extensi\'on a las ideas expresadas en
\cite{Cordova1}. En resumen, consideramos los conocimientos previos
de CALUNO, como sistema a analogar, sustentando fuertemente los
nuevos conocimientos y construcciones en CALVAV, como sistema o
escenario por conocer. Lo anterior, como ha sido mencionado, se
centra en los conceptos de derivada y diferenciabilidad.

Una de las fuentes de dificultad en el trabajo con funciones
bivariadas est\'a en las representaciones en tres dimensiones de sus
gr\'aficas y de los elementos (conceptos) que, inspirados en
\'estas, se pretendan destacar. Al  ser la noci\'on de dimensi\'on,
muy propia de las matem\'aticas superiores, de compleja definici\'on
y de usos muy diversos P\'AEZ \cite{Paez}, es claro que el uso de
alg\'un soporte visual para trabajar CALVAV puede ser un excelente
recurso pedag\'ogico complementario, facilitando la comprensi\'on y
el entendimiento (ANDRADE y MONTECINO \cite{Andrade}). Sin embargo,
nos limitamos en este art\'iculo a visualizar lo tridimensional en
el plano. Recordemos a PIAGET en \cite{Piaget}: {\it Como seres
humanos estamos limitados a visualizar la abstracci\'on de la
matem\'atica en dos dimensiones, que es diferente a nuestra
cotidianeidad}.

Las im\'agenes mentales son im\'agenes din\'amicas si en ellas los
objetos o algunos de sus elementos se desplazan GUTI\'ERREZ
\cite{Gutierrez}. En el CD, para presentar la derivada se recurre a
la imagen din\'amica de una recta secante sujeta a un punto en la
gr\'afica y donde el otro punto de intersecci\'on se va moviendo al
primero, de modo de crear la idea de estabilizaci\'on a una recta,
que despu\'es llamamos recta tangente. Las correspondientes
pendientes y su l\'imite, de existir, pasan a dar origen a la
derivada. Lo que se persigue con estas im\'agenes es la {\it
interpretaci\'on de informaci\'on figurativa}, es decir, la
comprensi\'on e interpretaci\'on de las  representaciones visuales
para extraer la informaci\'on que contienen, en contraste con el
{\it procesamiento visual} que es la conversi\'on de lo abstracto en
imagen visual, BISHOP \cite{Bishop}.

La principal novedad en el presente trabajo est\'a en el uso del
concepto de {\it pendiente de un plano} definido en \cite{Cordova2}
y la interpretaci\'on figurativa de la imagen de planos secantes (a
la gr\'afica de una superficie curva) sujetos a un punto com\'un,
estabiliz\'andose (din\'amica) para introducir el concepto de {\it
derivada de una funci\'on bivariada}. Finalmente, conjeturamos que
la existencia de esta derivada es equivalente a la cl\'asica
noci\'on de diferenciabilidad.

El plan del presente art\'iculo est\'a organizado como sigue. En la
Secci\'on 2, revisamos el concepto de analog\'ia y en particular el
de analog\'ia formal. En la Secci\'on 3, primera parte, resumimos
algunos elementos hist\'oricos y la interpretaci\'on de la imagen
din\'amica tras la construcci\'on de la derivada univariada para, en
una segunda, presentar lo que ser\'a su an\'aloga formal en varias
variables. En la Secci\'on 4, formulamos algunos ejemplos de
analog\'ias formales entre CALCUNO y CALVAV que se refieren a la
derivada. Cerramos con la Secci\'on 5, reservada a una discusi\'on y
ciertas conclusiones finales.

\section{Analog\'ias Formales}

Nuestra propuesta se enmarca en la ense\~nanza basada en
analog\'ias, dentro de la cual la {\it met\'afora} juega un rol
importante. Realizar una met\'afora es trasladar el significado
usual de una expresi\'on hacia otra en un contexto diferente, lo que
supone una relaci\'on de semejanza. Sin embargo, nos ocupar\'an las
{\it analog\'ias}, pues en ellas yace una creatividad cognitiva con
potencial de sistematizaci\'on en contraste con la met\'afora que
tiene un fin m\'as bien expresivo y por lo tanto sobrecargado
(P\'EREZ \cite{Perez}).

En la ense\~nanza de las ciencias exactas encontramos diversos
estudios que hacen buen uso de las analog\'ias (OLIVA \cite{Oliva},
DILBER \& DUZGUN \cite{Dilber},  FRIGO \& ADOLFO \cite{Frigo}) para
evidenciar la utilidad para los aprendizajes conceptuales y las
relaciones entre estos. Analog\'ias en contextos matem\'aticos
espec\'ificos son: Teor\'ia de grupos (SCHLIMM \cite{Schlimm}),
\'Algebra lineal (HOCHWALD \cite{Hochwald}), Ecuaciones
diferenciales (ABRAHAMSON \cite{Abrahamson}), entre otros.

Las analog\'ias tienen una etapa de metaforizaci\'on en la que se
encuentra la tarea de visualizaci\'on diferida, por medio de la cual
es posible atribuir caracter\'isticas a un dominio en aprendizaje
(denominado {\it Objetivo}), que en inicio no se manifiestan de
forma expl\'icita, sino que forman parte de un dominio ya conocido
(denominado {\it Base}), y son aplicados a este dominio Objetivo
(POCHULO \cite{Pochulu}). Cuando la Base est\'a fuera de las
matem\'aticas y el Objetivo dentro de \'estas, LAKOFF \& NU\~NEZ
\cite{Lakoff}) habla de met\'aforas (conceptuales) tipo {\it
Grounding}, traduscamos ``aterrizadas". Sin embargo, las que
presentaremos son del tipo {\it Linking} o ``enlazadas", en que
ambos dominios Base y Objetivo est\'an dentro de las matem\'aticas.

No obstante lo anterior, el uso de las analog\'ias en la ense\~nanza
posee tambi\'en sus limitaciones y no se recomienda su uso
asistem\'atico o improvisado pues como consecuencia puede dar lugar
a concepciones err\'oneas de conocimiento (DUIT \cite{Duit}).

En consideraci\'on con el cuidado que debe tenerse para introducir
nuevos conocimientos por medio de estas herramientas, GLYNN en
\cite{Glynn} desarrolla un modelo basado en la ense\~nanza con
analog\'ias (Teaching with Analogies, TWA de ahora en adelante) que
exige seguir una secuencia de etapas, de forma rigurosa y
sistem\'atica, con el objetivo de lograr un aprendizaje conceptual
exitoso en los estudiantes; como motor o coraz\'on de una
analog\'ia, una de las etapas apunta al mapeo que vincula An\'alogo
y T\'opico, instancia en la cual es posible apoyarse en la Teor\'ia
de mapeo estructural de la analog\'ia (SMTA por su denominaci\'on
original en ingl\'es, elaborada por GENTNER  \cite{Gentner83} y
GENTNER \& MARKMAN \cite{Gentner97}) que otorga simbolismo y
formalidad. Variante de dichos modelos es el Modelo Anal\'ogico
Formal expuesto en \cite{Cordova2}.

Relacionado a la ense\~nanza del c\'alculo diferencial (en
particular en varias variables) destaca el art\'iculo de G\'OMEZ \&
DELGADO \cite{Gomez} que aborda dicho contenido por medio de un
propuesta que guarda ribetes muy interesantes y vinculantes con el
presente estudio.

Por formalizaci\'on de una analog\'ia en \cite{Cordova2} se entiende
un modo de vincular e introducir alg\'un grado de rigor
matem\'atico, al menos en el lenguaje, en la definici\'on de las
asociaciones entre objetos del dominio Base con el dominio por
conocer u Objetivo.

Para definir una {\it analog\'ia formal} se requiere primero
reconocer dos conjuntos de objetos $X$ e $Y$, que pertenecen
respectivamente al dominio Base (denotado $S(X)$) y al Objetivo
(denotado $S(Y)$),  que tienen una asociaci\'on digamos ``intuitiva"
$T:X \to Y$.

Para analog\'ias de enlace, $S(X)$ y $S(Y)$ son dominios, \'areas o
teor\'ias de la matem\'atica y, sus respectivos conjuntos $X$ e $Y$,
son objetos matem\'aticos (e.g., conceptos o proposiciones) de
dichos dominios. Respecto al alcance, $S(Z)$ al menos es el conjunto
de proposiciones con sentido que involucran a los objetos de $Z\in
\{X,Y\}$. Adem\'as, la identificaci\'on $T:X \to Y$, se asume una
funci\'on inyectiva.

Entonces, una analog\'ia se denominar\'a analog\'ia formal cuando la
correspondencia entre proposiciones $p$ y $q$ de los respectivos
sistemas $S(X)$ y $S(Y)$ a la que refiere se puede expresar por
medio de una funci\'on proposicional $F$ en $n$-variables y un
arreglo de $n$-objetos $x_1, x_2,\cdots, x_n$ de $X$, tal que:
$$
p:\, F(x_1, x_2,\cdots, x_n)  \wedge q:\, F(T(x_1), T(x_2),\cdots,
T(x_n)),
$$
es verdadera.

Entonces, como se expresa en CORDOVA-LEPE et al. \cite{Cordova2}:
\begin{center}
{\it ``... obtener una analog\'ia es el acto de descubrir una
proposici\'on y su valor de verdad, un `consecuente'  $Q$ en un
dominio $S(Y)$, a trav\'es de la lectura de esta verdad en su
`antecedente', una proposici\'on $P$ del sistema $S(X)$. La
proposici\'on $P$ de $S(X)$ es para el alumno un saber (conceptual y
proposicional) ya conocido y cercano, luego mediante la extensi\'on
a una relaci\'on entre proposiciones (transducci\'on $P \to Q$)
desde una funci\'on de semejanza base $T: X \to Y$, que debe ser
expl\'icita,  \'el obtiene un conocimiento $Q$ en $S(Y)$".}
\end{center}

De particular inter\'es para nuestros fines son las analog\'ias que
se introducen entre (a) $S(X)$, la Geometr\'ia Anal\'itica de la
Recta en el Plano (GARP) y (b) $S(Y)$, la Geometr\'ia Anal\'itica
del Plano en el Espacio (GAPE). En que en $X$ están los conceptos de
{\it pendiente de una recta} e {\it intercepto de una recta}
asociados al gr\'afico de una funci\'on lineal af\'in $y=ax+b$.
Mientras que los objetos de $Y$  son los conceptos de {\it pendiente
de un plano} e {\it intercepto de un plano} asociados al gr\'afico
de una funci\'on $y=(a_{1},a_{2})\cdot (x_{1},x_{2})+b$. La
introducci\'on del concepto de pendiente de un plano se justifica
por semejanza sint\'actica con la notaci\'on $y=\vec{a}\cdot
\vec{x}+b$, donde $\vec{a}=(a_{1},a_{2})$ y $\vec{x}=(x_{1},x_{2})$,
para la funci\'on planar y tambi\'en por semejanza sem\'antica con
el rol geom\'etrico de $\vec{a}$ en la gr\'afica de dicha funci\'on.

\section{Secantes, Tangentes y Derivadas}
\subsection{Problema de la recta tangente en el plano cartesiano}

Con frecuencia los estudiantes reci\'en iniciados en los cursos de
c\'alculo tienen el preconcepto de recta tangente a una curva en uno
de sus puntos muy similar a la que los matem\'aticos griegos
cl\'asicos ten\'ian, esto es, como la recta que ``toca"  la curva en
dicho punto, pero sin cortarla. De hecho, para Euclides una recta
tangente a una circunferencia era la recta con un solo punto en
com\'un con dicha circunferencia.

Se les atribuye a Fermat y Descartes, siglo XVII, a trav\'es del
llamado {\it m\'etodo de las tangentes}, ser los primeros en abordar
el problema de encontrar, m\'as bien construir por medio de
algoritmos algebraicos, la recta tangente a una curva, ver ALARC\'ON
\& SUESC\'UN \cite{Alarcon1} y ALARC\'ON et al. \cite{Alarcon2}. Son
Newton y Leibniz los que introducen procedimientos ya m\'as
generales. Newton, asoci\'o la tangente con el vector velocidad de
una part\'icula, en cierto instante, al recorrer cierta curva, el
{\it m\'etodo cin\'etico}, ver SUZUKI \cite{Suzuki} y WOLFSON
\cite{Wolfson}. Algunos ven a aproximaci\'on de Leibniz en
conexi\'on con su metaf\'isica, principalmente con su``teor\'ia de
las monadas", ver MATEUS \cite{Mateus}.

Los matem\'atica tard\'o en adoptar la idea de recta tangente a una
curva como recta l\'imite de una familia de rectas secantes. Hubo
que esperar hasta Bolzano (1817) para tener indicios de la derivada
definida como un l\'imite. Seis a\~nos despu\'es, hacia 1823, fue
Cauchy quien  describi\'o la derivada, m\'as o menos como la
conocemos hoy, en  ``{\it Resum\'e des le\c{c}ons sur le calcul
infinitesimal}'', en t\'erminos del concepto de l\'imite.

En los procesos de aula, en cursos de c\'alculo en una variable, la
verbalizaci\'on y visualizaci\'on de los estudiantes del eslogan
\begin{itemize}
\item[{\bf p}:] {\it La recta tangente a una curva en un punto es el l\'imite de una familia de rectas secantes},
\end{itemize}
para una posterior matematizaci\'on, es algo que los profesores
desear\'iamos. En otras palabras la afirmaci\'on {\bf p}, con toda
su informalidad, resulta ser muy orientadora.

Nos gustar\'ia, en el contexto de c\'alculo en varias variables,
tener un ``gatillador" verbal y visual de potencial similar. En este
trabajo, queremos construir un an\'alogo a $P$, desde la
afirmaci\'on
\begin{itemize}
\item[{\bf q}:]
{\it El plano tangente a una superficie en un punto es el l\'imite
de una familia de planos secantes},
\end{itemize}
pues podr\'ia permitirnos una construcci\'on matem\'atica o
extensi\'on, natural al caso de varias variables con posibilidades
did\'acticas abiertas.

\vspace{5mm}
\noindent{\bf La derivada escalar}:
En lo que sigue de esta subsecci\'on, revisaremos la derivada
matem\'atica estandar tras la afirmaci\'on {\bf p}, con el fin de
extrapolar posteriormente una construcci\'on matem\'atica de
naturaleza semejante, pero ahora desde la proposici\'on {\bf q}.

Dada una funci\'on $f: I \to \mathbb{R}$, $I$ un intervalo abierto,
y un punto $x_{0}\in I$, un cambio en $x_{0}$ en $h$ unidades
determina el valor $x_{0}+h$. Notemos que $x_{0}+h$ est\'a en el
$I$, si $|\,h\,|$ es suficientemente peque\~no. El valor original
$x_{0}$ de la variable independiente y el nuevo $x_{0}+h$,
determinan respectivas im\'agenes, valores en la dependiente,
$f(x_{0})$ y $f(x_{0}+h)$. Al querer comparar las variaciones de las
variables, esto es, $h$ y $f(x_{0}+h)-f(x_{0})$, tenemos la
alternativa de cuantificar cuantas veces es la segunda por cada
unidad de la primera, mediante el cuociente
\begin{equation}
m(x_{0},h)=\frac{f(x_{0}+h)-f(x_{0})}{h}.
\end{equation}

Lo importante es que $m(x_{0},h)$ es interpretable geom\'etricamente
como la pendiente de la recta secante $L_{PQ}$ a la gr\'afica de $f$
en los puntos $P(x_{0},f(x_{0}))$ y $Q(x_{0}+h,f(x_{0}+h))$, ver
Fig. 1.

\begin{figure}[h!]
\begin{center}
\begin{tabular}{cc}
\includegraphics[height=4.1cm]{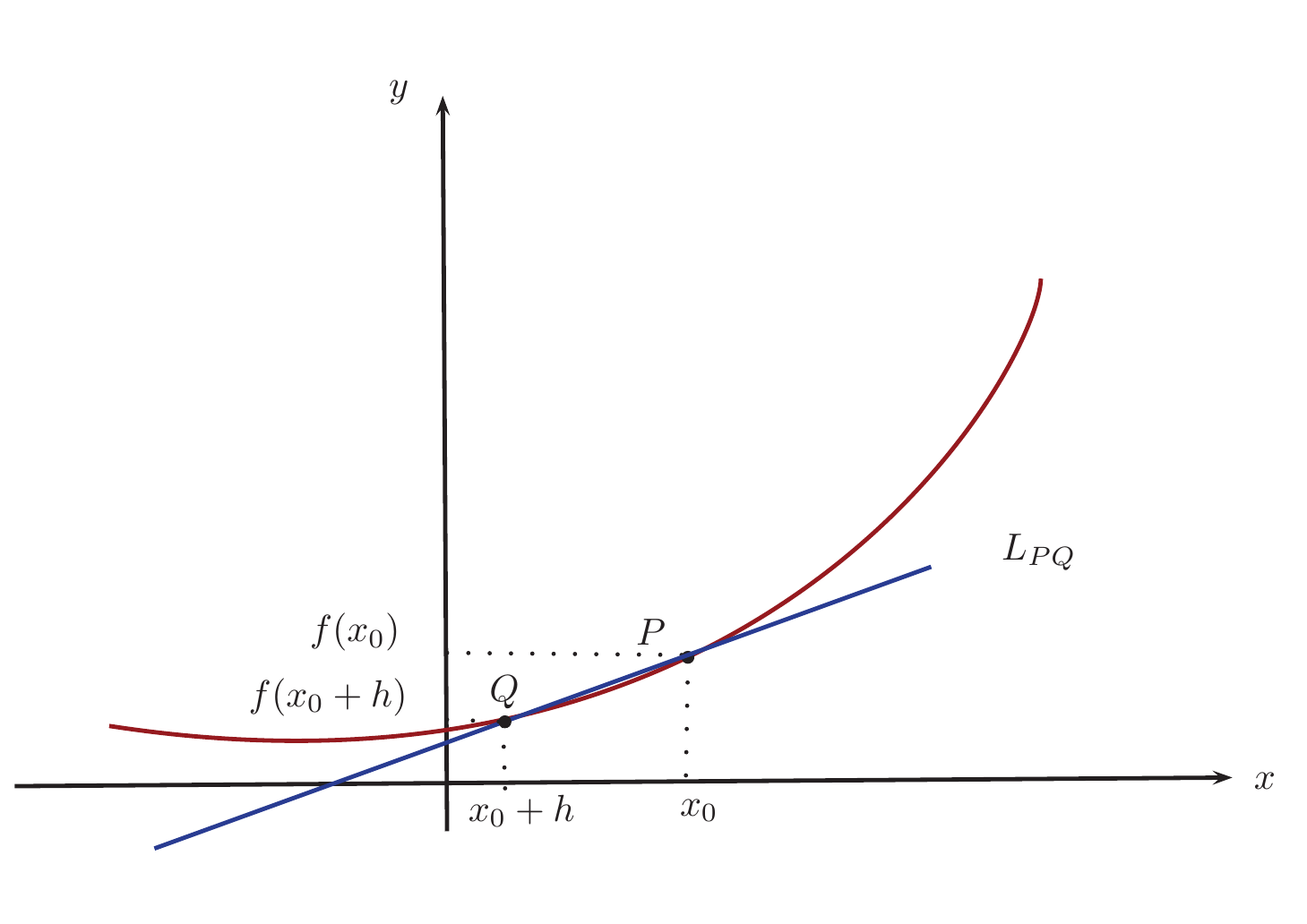} &
\includegraphics[height=4.2cm]{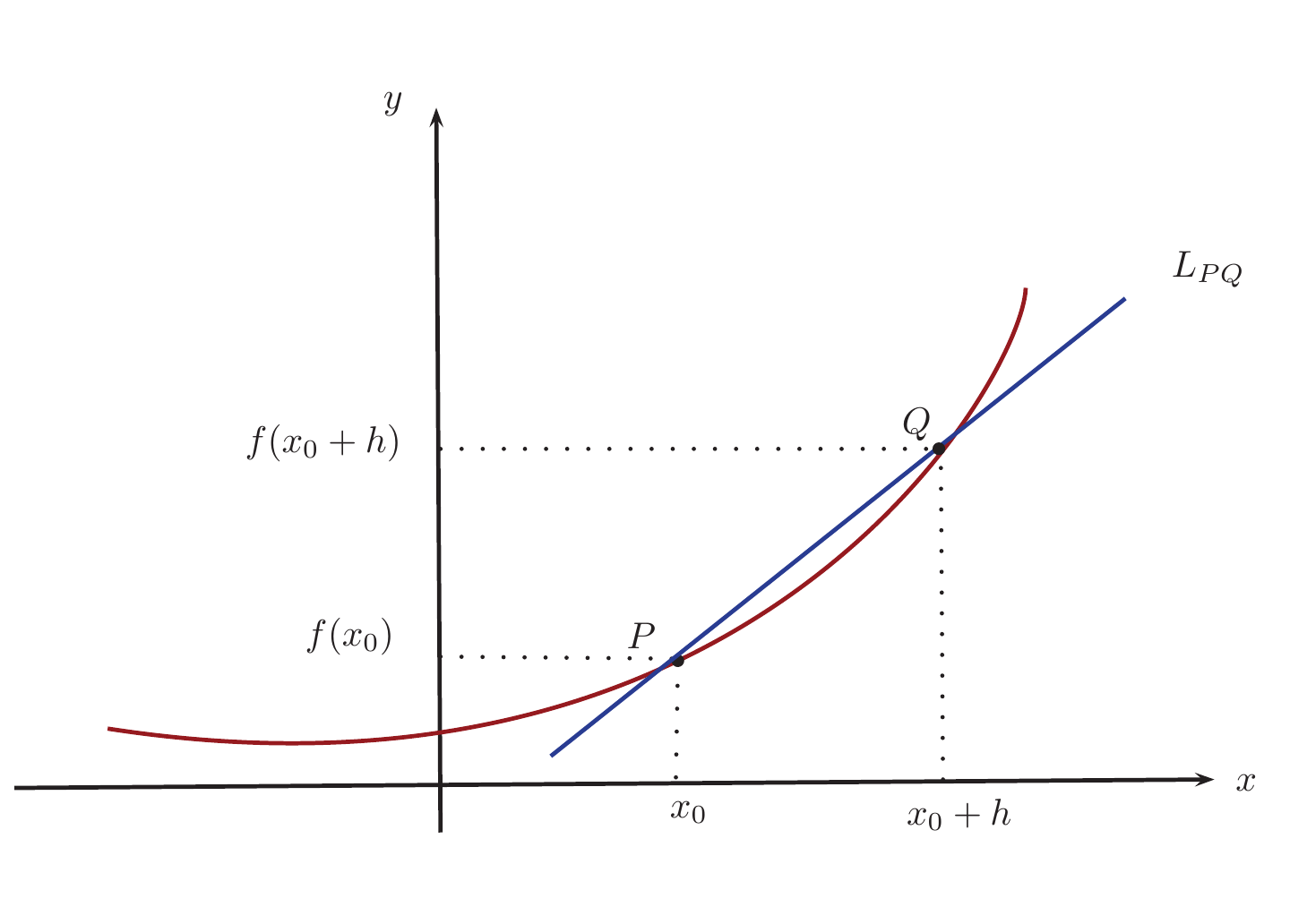}
\end{tabular}
\caption{Ilustración de la recta secante $L_{PQ}$ de pendiente
$m(x_{0},h)$ a la gráfica de $f$. Izquiera: $h<0$. Derecha: $h>0$.}
\end{center}
\end{figure}

El m\'etodo din\'amico se sustenta en lograr la visualizaci\'on de
lo que ocurre con la rectas $L_{PQ}$ (una im\'agen din\'amica)
cuando el punto $Q$ se mueve (tiende) hacia el punto $P$, el que
permanece fijo. Con diagramas de curvas suaves, se conf\'ia perfilar
en la imaginaci\'on la convergencia de las rectas $L_{PQ}$ a una
recta $L_{T}$ que recoge la idea de ``tocar" la gr\'afica de la
funci\'on en $P$ y, por lo que, merecer\'ia, por supuesto de
existir, el nombre de recta tangente, ver Fig. 2. Si las rectas
secantes se estabilizan en $L_T$, las pendientes de estas debieran
tambi\'en converger a lo que debiera ser la pendiente $m$ de $L_T$,
es decir,
$$
m(x_{0},h) \to m \,\,\,\mbox{cuando}\,\,\,h \to 0.
$$
Lo que corresponde, como ya sabemos, a la definici\'on de la {\it
derivada de la funci\'on} $f$ {\it en} $x_{0}$, si es que este
l\'imite existe, caso en que denotamos $m$ como $f'(x_{0})$.

\begin{figure}[h!]
\begin{center}
\includegraphics[height=4.2cm]{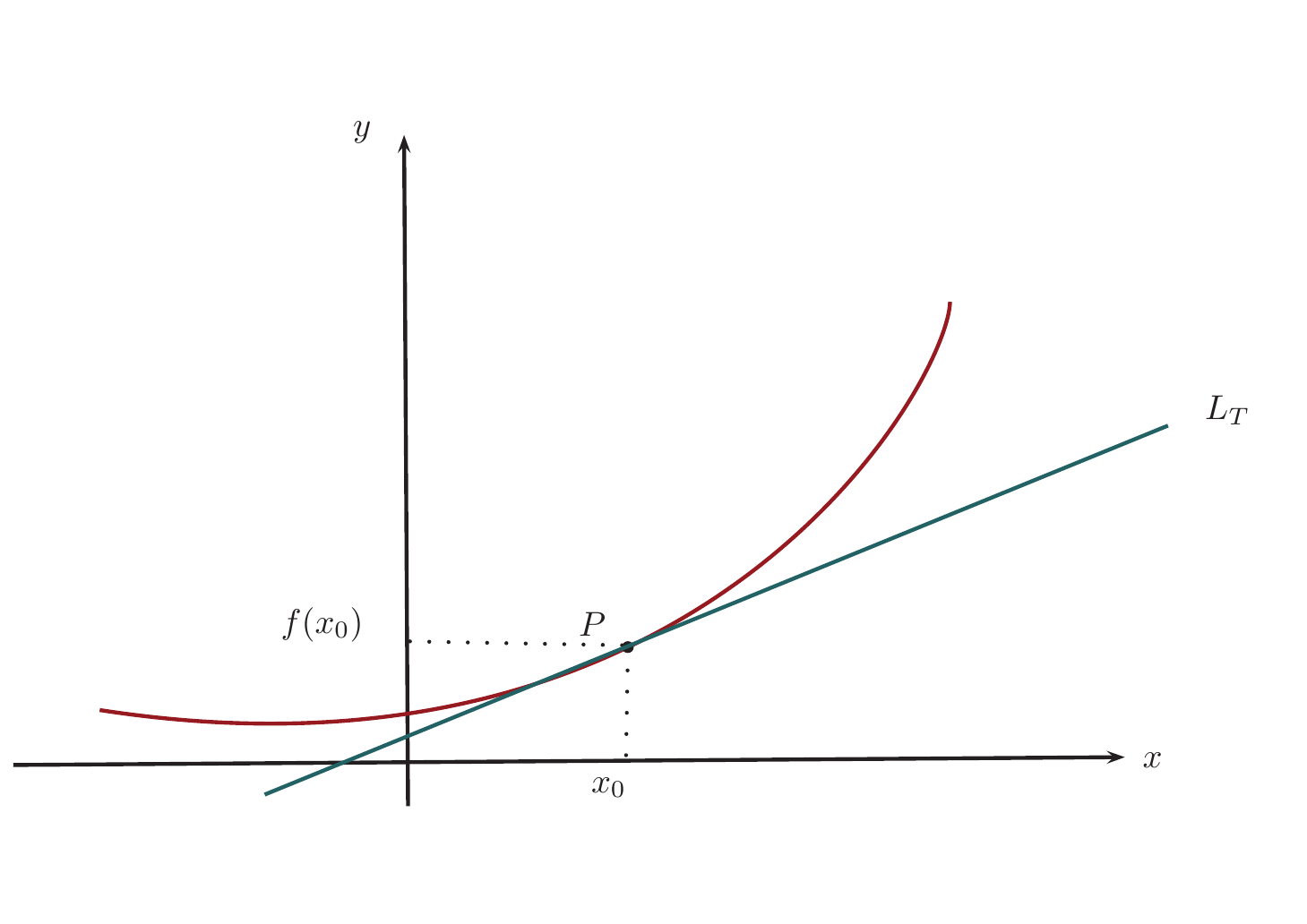}
\caption{Las rectas secantes $L_{PQ}$ se ``estabilizan" en una recta
$L_T$ que se denomina {\it recta tangente}.}
\end{center}
\end{figure}

\subsection{Problema del plano tangente en el espacio}

Manteniendo la esencia de la situaci\'on unidimensional, durante los
siglos XVIII y XIX se desarroll\'o el CD considerando ahora
funciones de dos y tres variables. Sin embargo, gran parte de la
terminolog\'ia empleada en las definiciones de l\'imites y
continuidad de funciones de dos o tres variables fueron introducidas
por el matem\'atico Karl Weierstrass, quien es considerado el {\it
padre} del an\'alisis moderno, \cite{Larson}.

Por otro lado, aunque Newton, Jean y Nicolaus Bernoulli hab\'ian
logrado la diferenciaci\'on de funciones de dos variables, fue con
los matem\'aticos Alexis Antoine de Bertis, Clairaut, d'Alembert y
Euler, donde esta teor\'ia alcanza su perfecci\'on. La idea era
simple: efectuar la derivada de una variable dejando el resto de
ellas como constantes \cite{Ruiz}. M\'as a\'un, entre $1730$ y
$1760$, Euler y d'Alembert publicaron de manera independiente una
gran cantidad de art\'iculos relacionados a problemas de equilibrio,
movimiento de flu\'idos y cuerdas vibrantes, utilizando una buena
parte de la teor\'ia asociada a derivadas parciales \cite{Larson}.
Estas teor\'ias se extendieron a la derivada direccional en el a\~no
de $1925$, gracias a las contribuciones de G\^ateaux \cite{Ramirez}.
Cabe destacar tambi\'en los desarrollos relacionados a la
formalizaci\'on del c\'alculo vectorial para la investigaci\'on en
F\'isica, culmin\'andose con el aporte de trabajos elaborados por
Hamilton, Grassmann, Tait y algunos investigadores de finales del
siglo XIX, donde aparecen nuevas nociones vectoriales, tales como
vector tangente y gradiente.

Por otra parte, el impacto de la Revoluci\'on Francesa en la ciencia
francesa y europea fue muy grande, concret\'andose en la reforma
educativa m\'as importante realizada desde el Renacimiento mediante
la creaci\'on de la {\em \'Ecole Normale Sup\'erieure, la \'Ecole de
M\'edecine} y la {\em \'Ecole Polytechnique}, donde uno de los
fundadores de esta \'ultima fue  Monge, el cual se suele
caracterizar como el primer especialista en geometr\'ia \cite{Ruiz},
en donde en su obra {\em G\'eometrie descriptive} estudia sobre
planos normales y tangentes a superficies curvas.

\vspace{5mm}
\noindent{\bf La derivada en varias variables}:

\begin{figure}[h!]
\begin{center}
\includegraphics[height=06cm]{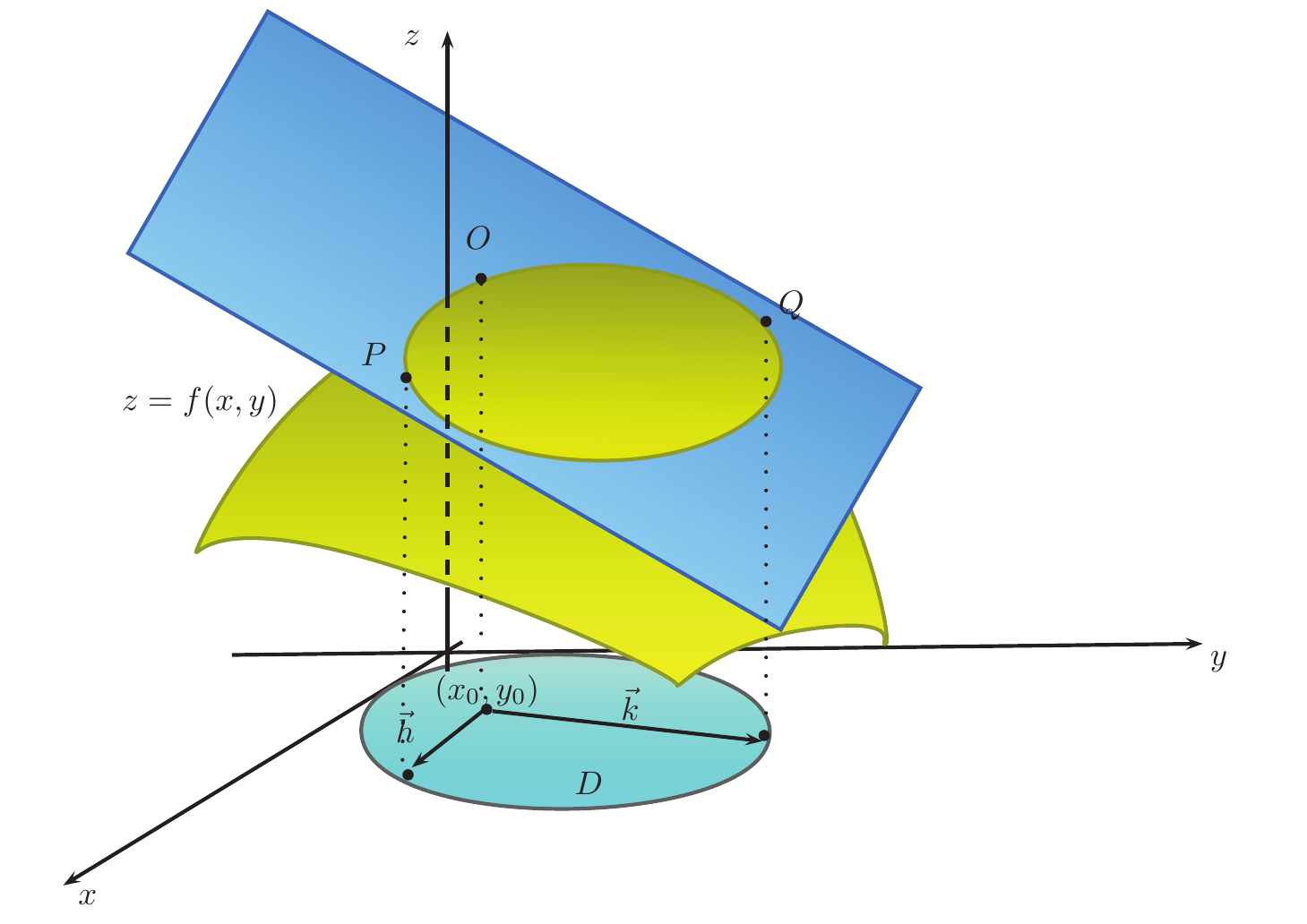}
\caption{El punto del dominio $(x_{0},y_{0})$ y los vectores
$\vec{h}$ y $\vec{k}$ determinan un plano secante a la gr\'afica de
$f$.}
\end{center}
\end{figure}

Consideremos ahora una funci\'on real de dos variables $z=f(x,y)$
sobre un cierto dominio abierto $D\subseteq \mathbb{R}^{2}$. Dado
$(x_{0},y_{0})$ en $D$, tenemos dos grados de libertad para variar
este punto, por lo que, definamos dos vectores linealmente
independientes
$$
\vec{h}=(h_{1},h_{2}) \qquad \mbox{ y } \qquad \vec{k}=(k_{1},k_{2})
$$
con $||\,\vec{h}\,||$ y $||\, \vec{k}\,||$ suficientemente pequeños
para que $(x_{0},y_{0})+\vec{h}$ y $(x_{0},y_{0})+\vec{k}$
pertenezcan al dominio $D$.

Tenemos las im\'agenes $f(x_{0},y_{0})$, $f((x_{0},y_{0})+\vec{h})$
y $f((x_{0},y_{0})+\vec{k})$, los que determinan puntos $$
O((x_{0},y_{0}),f(x_{0},y_{0})), \,\,\,
P((x_{0},y_{0})+\vec{h},f((x_{0},y_{0})+ \vec{h})) \,\,\, \mbox{y}
\,\,\, Q((x_{0},y_{0})+\vec{k},f((x_{0},y_{0})+ \vec{k}))
$$
en la gr\'afica de $f$. Notemos que $O$, $P$ y $Q$ definen un plano
secante a la gr\'afica de $f$, ver Fig. 3, de ecuaci\'on
(punto-pendiente):
$$
z=(a,b) \cdot [(x,y)-(x_{0},y_{0})]+f(x_{0},y_{0}),
$$
donde el par $(a,b)$ es lo que mide y se denomina  {\em{pendiente}}
de dicho plano y que también denotaremos por
$m[(x_{0},y_{0}),\vec{h}, \vec{k}\,]$.

Dado lo anterior, al remplazar en la ecuaci\'on de este plano los
puntos $P$ y $Q$, tenemos el siguiente sistema de ecuaciones para
las inc\'ognitas $a$ y $b$:
\begin{equation}
\left\{
\begin{array}{lll}
(a,b)\cdot \vec{h}&=a\,h_{1}+b\,h_{2}=f((x_{0},y_{0})+\vec{h}\,)-f(x_{0},y_{0})=&\Delta_{\vec{h}} f(x_{0},y_{0})\\
(a,b)\cdot
\vec{k}&=a\,k_{1}+b\,k_{2}=f((x_{0},y_{0})+\vec{k}\,)-f(x_{0},y_{0})=&\Delta_{\vec{k}}
f(x_{0},y_{0})
\end{array}
\right.
\end{equation}

Usando la Regra de Cramer obtenemos
\begin{eqnarray}
m[(x_{0},y_{0}),\vec{h}, \vec{k}\,]= \left(\frac{k_{2}\,
\Delta_{\vec{h}}f(x_{0},y_{0})-h_{2}\,\Delta_{\vec{k}}f(x_{0},y_{0})}{h_{1}
k_{2}-k_{1} h_{2}},\frac{h_{1}
\Delta_{\vec{k}}f(x_{0},y_{0})-k_{1}\Delta_{\vec{h}}f(x_{0},y_{0})
}{h_{1} k_{2}-k_{1} h_{2}}\right). \label{pendiente}
\end{eqnarray}

Cuando $\vec{h}$ y $\vec{k}$ var\'ian y tanto $||\,\vec{h}\,||$ como
$||\,\vec{k}\,||$ se hacen peque\~nos, en la Fig. 4 es posible
visualizar los planos secantes estabiliz\'andose a lo que
pordr\'iamos llamar el plano tangente a la gr\'afica de $f$ en el
punto $O$, su pendiente, que denotaremos $f'(x_{0},y_{0})$ se
obtiene, de existir, como
\begin{equation}
\label{definition} f'(x_{0},y_{0})=\lim_{(\,\vec{h},\,\vec{k}\,) \to
(\,\vec{0},\,\vec{0}\,)} m[(x_{0},y_{0}),\vec{h},\vec{k}].
\end{equation}
De existir $f'(x_{0},y_{0})$ se llamará la {\it derivada de $f$ en
$(x_{0},y_{0})$} o diremos que {\it $f$ es derivable} en tal punto.

\begin{figure}[h!]
\begin{center}
\includegraphics[height=06cm]{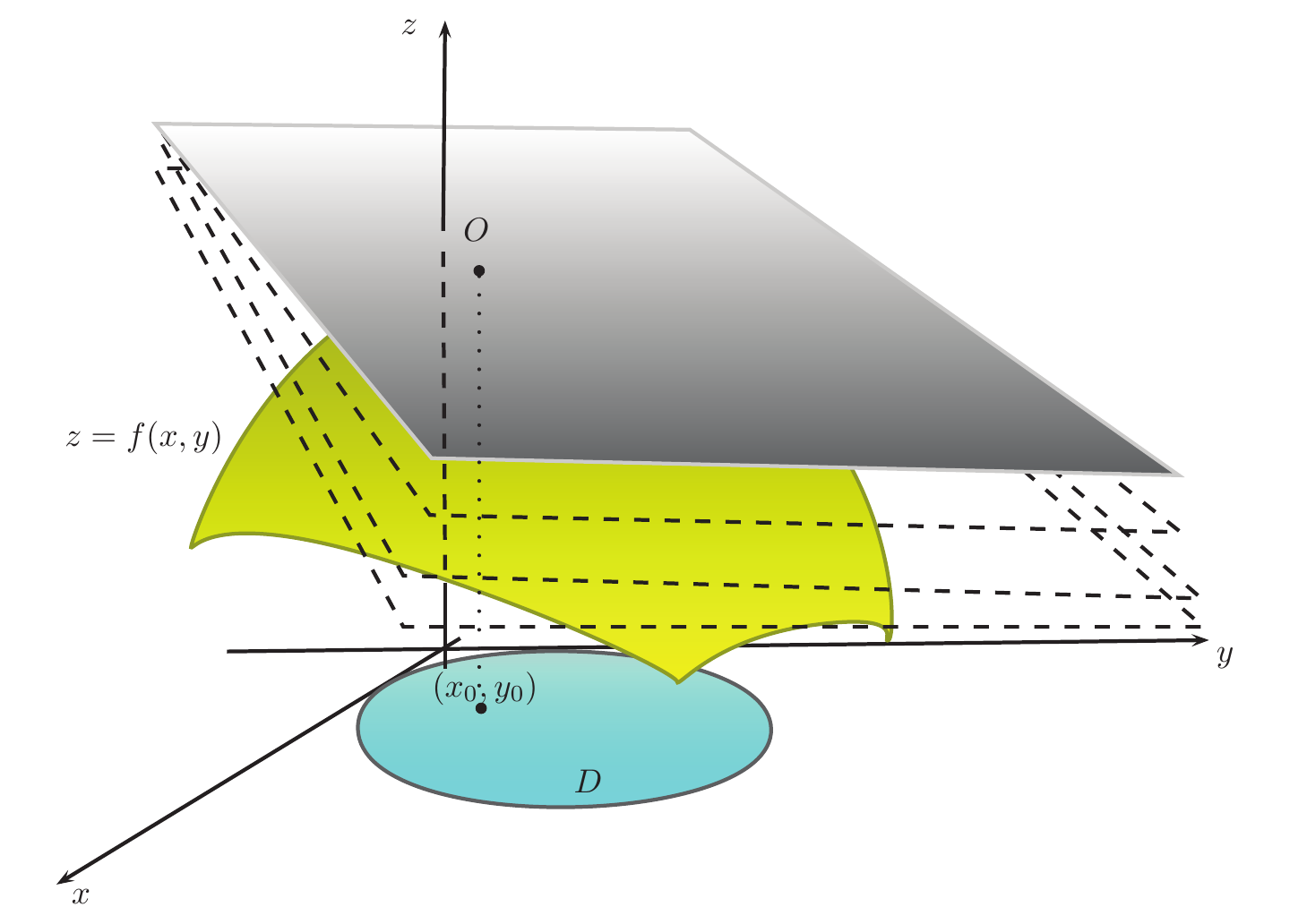}
\caption{Imagen que intenta representar planos secantes
``estabiliz\'andose" a lo que se denomina {\it plano tangente}.}
\end{center}
\end{figure}

\vspace{5mm}
\noindent{\bf Algunas propiedades de c\'alculo}:
Sean $f:D_{f} \subseteq \mathbb{R}^{2} \rightarrow \mathbb{R}$ y
$g:D_{g} \subseteq \mathbb{R}^{2} \rightarrow \mathbb{R}$ dos
funciones y supongamos que tenemos un par $(x_{0},y_{0}) \in D_{f}
\cap D_{g} $ donde $f'(x_{0},y_{0})$ y $g'(x_{0},y_{0})$ existen,
entonces se verifica que:
\begin{itemize}
\item[(a)] Si $f(x,y)=\alpha x + \beta y + \gamma$, entonces
$f'(x_{0},y_{0})=(\alpha,\beta)$ para todo trío de números reales
$\alpha$, $\beta$ y $\gamma$. En particular $x'=(1,0)$ e $y'=(0,1)$.

\item[(b)] Como $\Delta_{\vec{v}}$ es un operador lineal, se tiene la derivada tambi\'en, i.e., $(\alpha f+\beta g)'(x_{0},y_{0})=\alpha f'(x_{0},y_{0})+\beta g'(x_{0},y_{0})$, para todo par de n\'umeros reales $\alpha$ y $\beta$.
\end{itemize}

Con un poco m\'as de esfuerzo y recurriendo a los m\'etodos
similares de CALUNO en (\ref{pendiente}), se tiene:
\begin{itemize}
\item[(c)] La regla para la derivada de un producto es
$(fg)'(\cdot)=f'(\cdot)\,g(\cdot)+f(\cdot)\,g'(\cdot)$. Por ejemplo,
$(x^{2}y^{3})'$ es igual a $2xx'y^{3}+x^{2}3y^{2}y'$. Entonces por
(a) se obtiene $2xy^{3}(1,0)+x^{2}3y^{2}(0,1)$, esto es,
$(x^{2}y^{3})'=(2xy^{3},3x^{2}y^{2})$.

\item[(d)] La regla para la derivada de un cuociente es
$(1/f)'(\cdot)=-f'(\cdot)/f^{2}(\cdot)$. De modo que $(x/y)'$ es
igual a $x'(1/y)+x(-y'/y^{2})$. Es decir, $(x/y)'=(1/y,-x/y^{2})$.
\end{itemize}

\vspace{5mm}
\noindent{\bf Derivabilidad de $f$ implica $f'(\cdot)=\nabla
f(\cdot)$ }:
Suponemos que $f:D\to \mathbb{R}$ es una funci\'on derivable en el
sentido reci\'en introducido en $(x_{0},y_{0})\in D$. Entonces, en
(\ref{definition}) podemos condicionar los vectores linealmente
independientes $\vec{h}$ y $\vec{k}$ para el c\'alculo del l\'imite.

Notemos que si $\vec{h}=s\,(H_{1},H_{2})$ y
$\vec{k}=t\,(K_{1},K_{2})$ con
$||(H_{1},H_{2})||=||(K_{1},K_{2})||=1$, entonces
$m[(x_{0},y_{0}),\vec{h},\vec{k}]$ es igual a
$$
\frac{1}{H_{1}K_{2}-K_{1}H_{2}} \left( K_{2}\,
\frac{\Delta_{\vec{h}}f(x_{0},y_{0})}{s}- H_{2}\,
\frac{\Delta_{\vec{k}}f(x_{0},y_{0})}{t},\, H_{1}\,
\frac{\Delta_{\vec{k}}f(x_{0},y_{0})}{t}- K_{1}\,
\frac{\Delta_{\vec{h}}f(x_{0},y_{0})}{s} \right).
$$

Notemos que
$$
\lim_{s\to 0}\frac{\Delta_{\vec{h}}f(x_{0},y_{0})}{s} \quad \mbox{y}
\quad \lim_{t\to 0}\frac{\Delta_{\vec{k}}f(x_{0},y_{0})}{t},
$$
corresponden a las derivadas direccionales de $f$ en $(x_{0},y_{0})$
en las direcciones de $\vec{H}=(H_{1},H_{2})$ y
$\vec{K}=(K_{1},K_{2})$ respectivamente, esto es,
$$
\nabla f(x_{0},y_{0}) \cdot (H_{1},H_{2}) \quad \mbox{y} \quad
\nabla f(x_{0},y_{0}) \cdot (K_{1},K_{2}).
$$
Finalmente, tomando límite en la última expresión para
$m[(x_{0},y_{0}),\vec{h}, \vec{k}]$, obtenemos
$$
f'(x_{0},y_{0})=\lim_{(s,t)\to (0,0)} m[(x_{0},y_{0}),\vec{h},
\vec{k}]=\nabla f(x_{0},y_{0}).
$$

\vspace{5mm}

\noindent{\bf Ejemplo de funci\'on no derivable}:
Consideremos la funci\'on $f:D \subset \mathbb{R}^{2}\to \mathbb{R}$
definida como $f(x,y)=x^{2}y/(x^{4}+y^{2})$ si $(x,y)\neq (0,0)$ y
tal que $f(0,0)=0$.

Asumiremos que $\vec{h}=(h_{1},h_{2})=\lambda(1,\alpha)$ y
$\vec{k}=(k_{1},k_{2})=\lambda(1,\beta)$ con $\alpha \neq \beta$.
Entonces al variar $(0,0)$ en $\vec{h}$ o en $\vec{k}$ tenemos
$$
\Delta_{\vec{v}}f(0,0)=f(v_{1},v_{2})=\frac{v_{1}^{2}v_{2}}{v_{1}^{4}+v_{2}^{2}},\,\,\vec{v}\in
\{\vec{h},\vec{k}\,\}.
$$
Luego, $\Delta_{\vec{h}}f(0,0)=\lambda
\alpha/\{\lambda^{2}+\alpha^{2}\}$  y $
\Delta_{\vec{k}}f(0,0)=\lambda \beta/\{\lambda^{2}+\beta^{2}\}$. Al
remplazar en (\ref{pendiente}) tenemos
$$
m[(0,0),\vec{h}, \vec{k}]= \left(\frac{\alpha
\beta(\alpha+\beta)}{(\lambda^{2}+\alpha^{2})(\lambda^{2}
+\beta^{2})}   , \frac{\lambda^{2}-\alpha
\beta}{(\lambda^{2}+\alpha^{2})(\lambda^{2} +\beta^{2})} \right).
$$
Con lo que
$$
=\lim_{\lambda \to 0} m[(x_{0},y_{0}),\vec{h}, \vec{k}]=
\left(\frac{\alpha+\beta}{\alpha \beta},-\frac{1}{\alpha \beta}
\right),
$$
por lo tanto, el l\'imite no existe con independencia de $\vec{h}$ y
$\vec{k}$. De modo que la funci\'on no es derivable en $(0,0)$.

\section{Analog\'ias formales entre derivadas}
En esta secci\'on, usaremos el concepto introducido en
\cite{Cordova2} de {\it analog\'ia formal}, la  cual es una
especificaci\'on mediante el lenguaje matem\'atico, para una
noci\'on de modelo anal\'ogico cercana a la de GENTNER
\cite{Gentner89}). En mencionado trabajo, se ejemplifica con
analog\'ias entre la Geometr\'ia Anal\'itica de la Recta en el Plano
(GARP) y la Geometr\'ia Anal\'itica del Plano en el Espacio (GAPE).
Mediante la construcci\'on de los conceptos de {\it pendiente} e {\it intercepto} de un plano, se afirma que en GAPE \\
``{\it Dos \underline{planos} son paralelos si y s\'olo si sus \underline{pendientes} son iguales y no as\'i sus \underline{interceptos}}"\\
como an\'aloga a la proposici\'on base \\
``{\it Dos \underline{rectas} son paralelas si y s\'olo si sus
\underline{pendientes} son iguales y no as\'i sus
\underline{interceptos}}.

Notemos que en la secci\'on anterior, rescatamos la introducci\'on
en \cite{Cordova2} del concepto de pendiente de un plano para ligar
con la construcci\'on de la derivada para funciones bivariadas.

Introducimos los siguientes sistemas de conocimientos:
\begin{itemize}
\item $S(X)$: C\'alculo Diferencial en Una Variable (CALUNO) y
\item $S(Y)$: C\'alculo Diferencial en Varias Variables (CALVAV),
\end{itemize}
donde $X$ e $Y$ son respectivamente algunos objetos matem\'aticos de
CALUNO y CALVAV asociados mediante una funci\'on base $T:X \to Y$
uno a uno.

\vspace{5mm}
\noindent{\bf Ejemplo 1}:\\
Consideremos los conjuntos $X_{1}=\{\mbox{derivada},
\mbox{pendientes}, \mbox{rectas secantes}\}$ e
$Y_{1}=\{\mbox{gradiente}, \mbox{pendientes}, \mbox{planos
secantes}\}$. La funci\'on $T_{1}:X_{1} \to Y_{1}$ definida por
$T_{1}(\mbox{derivada})=\mbox{gradiente}$,
$T_{1}(\mbox{pendiente})=\mbox{pendiente}$ y $T_{1}(\mbox{recta
secante})= \mbox{plano secante}$, permite definir desde la funci\'on
l\'ogica
$$
\begin{array}{rl}
F_{1}(x_{1},x_{2},x_{3}):& x_{1}\,\,\mbox{es el l\'imite de las} \,\,x_{2}\,\,\mbox{de una familia}\\
& \mbox{de}\,\,x_{3}\,\,\mbox{con un punto en com\'un.}
 \end{array}
$$
Desde esta funci\'on determinamos las proposiciones an\'alogas
$P_{1}=F(x_{1},x_{2},x_{3})$ y
$Q_{1}=F(T_{1}(x_{1}),T_{1}(x_{2}),T_{1}(x_{3}))$ con
$x_{1}=\mbox{derivada}$, $x_{2}=\mbox{pendientes}$ y
$x_{3}=\mbox{rectas secantes}$. Esto es,
\begin{itemize}
\item $P_{1}$: {\it La derivada de una funci\'on es el l\'imite de las pendientes de una familia de rectas secantes con un punto en com\'un}.

\item
$Q_{1}$: {\it El gradiente de una funci\'on es el l\'imite de las
pendientes de una familia de planos secantes con un punto en
com\'un}.
\end{itemize}

\vspace{5mm}
\noindent{\bf Ejemplo 2}:\\
Consideremos ahora los conjuntos $X_{2}=\{\mbox{derivada},
\mbox{pendientes}, \mbox{recta tangente}\}$ e
$Y_{2}=\{\mbox{gradiente}, \mbox{pendientes}, \mbox{plano
tangente}\}$. La función $T_{2}:X_{2} \to Y_{2}$ definida por
$T_{2}(\mbox{derivada})=\mbox{gradiente}$,
$T_{2}(\mbox{pendiente})=\mbox{pendiente}$ y $T_{1}(\mbox{recta
tangente})= \mbox{plano tangente}$, permite definir desde la
funci\'on l\'ogica
$$
\begin{array}{rl}
F_{2}(x_{1},x_{2},x_{3}):& x_{1}\,\,\mbox{de una función define la}
\,\,x_{2}\,\,\mbox{de la }\,\,x_{3}.
 \end{array}
$$
Desde esta funci\'on determinamos las proposiciones an\'alogas
$P_{2}=F(x_{1},x_{2},x_{3})$ y
$Q_{2}=F(T_{1}(x_{1}),T_{1}(x_{2}),T_{1}(x_{3}))$ con
$x_{1}=\mbox{derivada}$, $x_{2}=\mbox{pendientes}$ y
$x_{3}=\mbox{rectas tangente}$. Esto es,
\begin{itemize}
\item $P_{2}$: {\it La derivada de una funci\'on define a la pendiente de la recta tangente}.

\item
$Q_{2}$: {\it El gradiente de una funci\'on define a la pendiente
del plano tangente}.
\end{itemize}

\vspace{5mm}

\noindent{\bf Otras posibilidades}:\\
Las definiciones tambi\'en pueden tener una construcci\'on a base de
analog\'ias formales. Es el caso cuando al hacer la identificaci\'on
curva con superficie y dimensi\'on dos con tres, tenemos las
siguientes lecturas: {\it Una recta (plano) que intersecta en al
menos dos (tres) puntos distintos a una curva se llamar\'a recta
(plano) secante}. Sin embargo, sabemos que la segunda lectura
requiere la \underline{no colinealidad} de los tres puntos, cosa que
por estructura no podemos (al menos no de una manera obvia) asimilar
como an\'alogo a un concepto del dominio base. La delimitaci\'on de
los alcances de una analog\'ia es un tema importante, respecto a
nuestra analog\'ia, que ahora podemos llamar entre derivadas de
CALUNO y CALVAV, tambi\'en ser\'a necesario establecerlas, pero no
es el momento en el presente trabajo.


\section{Discusi\'on}

El camino que se inici\'o en C\'ORDOVA-LEPE et al. \cite{Cordova2}
con la introduci\'on de la {\it analog\'ia formal}, que permiti\'o
definir qu\'e se entiende por {\it pendiente de un plano}, mediante
analog\'ia con la pendiente de una recta, en el presente trabajo
desemboc\'o en la construcci\'on del concepto de {\it derivada} para
una funci\'on bivariada. Proceso que es claramente generalizable a
funciones de dominios reales en dimensiones mayores. Nuestra
intuici\'on y primeras pruebas indica que la noci\'on de
derivabilidad introducida absorbe o es equivalente a la cl\'asica
diferenciabilidad para funciones de varias variables. Sin dudas se
abre un espacio de exploraci\'on matem\'atica.

Sin embargo, la conjetura de que esta derivabilidad ofrece una
alternativa did\'actica en comparaci\'on a lo que vemos en los
textos de estudio y realizamos en las aulas de enseñanza superior,
est\'a por ser validada. La invitaci\'on es a responder:
\begin{itemize}
\item ¿ Qu\'e diferencias en la comprensi\'on conceptual es posible establecer entre un estudiante que ha abordado el tema de la derivaci\'on de una
funci\'on escalar de dos variables a trav\'es de un modelo
anal\'ogico formal, y un estudiante que lo ha abordado por medio de
una ense\~nanza tradicional?

\item ¿Qu\'e justificaci\'on entrega un estudiante, que ha abordado el
concepto de gradiente por medio de un modelo anal\'ogico formal, al
resolver una situación problemática de optimizaci\'on, y qu\'e
diferencias se evidencian en contraste con los argumentos emitidos
por un estudiante que ha abordado el contenido de forma tradicional?
\end{itemize}

Puede resultar orientador en el trabajo de aula respetar la etapas
jer\'arquicas que GLYNN \cite{Glynn} clasifica para la ense\~nanza
con analog\'ias: (1) Introducir el concepto objetivo (la derivada de
funciones bivariadas), (2) Realizar un diagnóstico de lo que
recuerdan los estudiantes del concepto análogo (la derivada de
CALUNO), (3) Identificar las características relevantes de los
conceptos objetivo y análogo (l\'imite de las pendientes de una
familia de rectas, respectivamente planos, secantes), (4)   Conectar
las caracter\'isticas similares entre ambos conceptos (recta con
plano, pendiente de recta con pendiente de plano, recta secante con
plano secante), (5) Indicar las limitaciones de la analogía (por
ejemplo, la independencia lineal de las direcciones en el l\'imite)
y (6) Establecer conclusiones en relación al concepto objetivo (las
propiedades de la derivabilidad y su relaci\'on con el concepto de
diferenciabilidad). Se trata de desarrollar los conceptos
involucrados, nuevos y basales, de una forma tal, que estos les
resulten significativos a los estudiantes.

Que los profesores o autores de texto utilicen las t\'ecnicas de
comparaci\'on con los conocimientos previos para introducir nuevas
ideas, es recurrente. Aunque, a veces utilizan las analog\'ias en
forma automatizada y sin estar conscientes de su uso. Es mejor, para
no perjudicar el aprendizaje de los estudiantes, por ejemplo
desarrollando concepciones err\'oneas (DUIT \cite{Duit}), que usen,
las analogías o las analog\'ias formales de manera m\'as
sistem\'atica. Las investigaciones deben avanzar en cuanto a
dise\~nar y validar modelos estrat\'egicos y sistem\'aticos para
ense\~nar por medio de analogías (GLYNN \cite{Glynn}).

\vspace{5mm}

\noindent{\bf Agradecimientos}: Los dos primeros autores agradecen
el patrocinio y auspicio del Plan de Mejoramiento Institucional
UCM1310, MINEDUC, Chile.

\end{document}